\documentstyle[11pt]{article}
\title{\LARGE  On the best constant of    Hardy--Sobolev Inequalities}

\author{\Large Adimurthi$^{1}$ \& Stathis Filippas$^{2,4}$
 \&  Achilles Tertikas$^{3,4}$  \\
                                                      \\
     TIFR center P.O. Box 1234$^{1}$,  \\
         Bangalore 560012, India \\
         aditi@math.tifrbng.res.in \\
  \\
        Department of Applied Mathematics$^{2}$ \\
         University of Crete,
        71409 Heraklion,  Greece \\
        filippas@tem.uoc.gr\\
                                          \\
 Department of Mathematics$^{3}$ \\
         University of Crete,
         71409 Heraklion,  Greece \\
          tertikas@math.uoc.gr\\
                                      \\
        Institute of Applied and Computational Mathematics$^4$, \\
        FORTH, 71110 Heraklion, Greece \\
    \\ }

\begin{document}
\date{}
\newcommand{\ana}{\nabla}
\newcommand{\R}{I \!  \! R}
\newcommand{\HH}{I \!  \! H}
\maketitle

\newcommand{\be}{\begin{equation}}
\newcommand{\ee}{\end{equation}}
\newcommand{\bea}{\begin{eqnarray}}
\newcommand{\eea}{\end{eqnarray}}
\newcommand{\la}{\label}
\newcommand{\xa}{\alpha}
\newcommand{\xb}{\beta}
\newcommand{\xg}{\gamma}
\newcommand{\xG}{\Gamma}
\newcommand{\xd}{\delta}
\newcommand{\xD}{\Delta}
\newcommand{\xe}{\varepsilon}
\newcommand{\xz}{\zeta}
\newcommand{\xh}{\eta}
\newcommand{\th}{\theta}
\newcommand{\Th}{\Theta}
\newcommand{\xk}{\kappa}
\newcommand{\xl}{\lambda}
\newcommand{\xL}{\Lambda}
\newcommand{\xm}{\mu}
\newcommand{\xn}{\nu}
\newcommand{\ks}{\xi}
\newcommand{\KS}{\Xi}
\newcommand{\xp}{\pi}
\newcommand{\xP}{\Pi}
\newcommand{\xr}{\rho}
\newcommand{\xs}{\sigma}
\newcommand{\xS}{\Sigma}
\newcommand{\xf}{\phi}
\newcommand{\xF}{\Phi}
\newcommand{\ps}{\psi}
\newcommand{\PS}{\Psi}
\newcommand{\xo}{\omega}
\newcommand{\xO}{\Omega}
\newcommand{\Ren}{ I \! \! R^N}
\newcommand{\Real}{ I \! \! R}
\newcommand{\spc}{{\cal C}_0^{\infty} (\xO)}
\newcommand{\sph}{H_0^1 (\xO)}
\newcommand{\spw}{W_0^{1,2} (B_1 ;\;  \weight )}
\newcommand{\iuu}{\int_{\Ren} |\nabla u|^2\,dx}
\newcommand{\ioo}{\int_{\Ren} Vu^2\,dx}
\newcommand{\Br}{B_r}
\newcommand{\bBr}{\partial \! B_r}
\newcommand{\ra}{\rightarrow}
\newcommand{\rft}{\rightarrow +\infty}
\newcommand{\bin}{\int_{\bBr}}
\newcommand{\weight}{|x|^{-(n-2)}}
\newcounter{newsection}
\newtheorem{theorem}{Theorem}[section]
\newtheorem{lemma}[theorem]{Lemma}
\newtheorem{prop}[theorem]{Proposition}
\newtheorem{coro}[theorem]{Corollary}
\newtheorem{defin}[theorem]{Definition}
\newcounter{newsec} \renewcommand{\theequation}{\thesection.\arabic{equation}}

\begin{abstract}
We obtain  the sharp  constant for the   Hardy-Sobolev inequality involving the distance to the origin.
This inequality is equivalent to a limiting Caffarelli--Kohn --Nirenberg inequality. In  three dimensions,
in certain  cases the sharp constant coincides with the best Sobolev constant.
\end{abstract}


\setcounter{equation}{0}
\section{Introduction}

The standard  Hardy inequality  involving the distance to the origin,
 asserts that  when  $n \geq 3$  and     $u \in C^{\infty}_{0}(\R^n)$
one has \be\la{0.1} \int_{\R^n} |\nabla u|^2 dx \geq
\left(\frac{n-2}{2} \right)^2   \int_{\R^n} \frac{u^2}{|x|^2} dx.
\ee The constant  $\left(\frac{n-2}{2} \right)^2$ is the best
possible  and remains the same if we replace  $u \in
C^{\infty}_{0}(\R^n)$  by  $u \in C^{\infty}_{0}(B_1)$, where $B_1
\subset \R^n$ is the unit ball centered at zero.  Brezis and
V\'{a}zquez  \cite{BV} have improved it by establishing that for
$u \in C^{\infty}_{0}(B_1)$, \be\la{0.2} \int_{B_1} |\nabla u|^2
dx \geq   \left(\frac{n-2}{2} \right)^2   \int_{B_1}
\frac{u^2}{|x|^2} dx + \mu_1  \int_{B_1} u^2 dx, \ee where  $\mu_1
= 5.783...$ is the first eigenvalue of the Dirichlet Laplacian of
the unit disc in $\R^2$. We note that  $\mu_1$   is the best
constant  in (\ref{0.2})
  independently of the dimension $n \geq 3$.

When taking distance to the boundary,  the following Hardy inequality  with best constant is also well known
for $n \geq 2$ and  $u \in C^{\infty}_{0}(B_1)$,
\be\la{0.3}
\int_{B_1} |\nabla u|^2 dx \geq \frac14   \int_{B_1} \frac{u^2}{(1-|x|)^2} dx.
\ee
Similarly to (\ref{0.2}) this has  also been improved  by Brezis and Marcus in  \cite{BM}  by
proving that
\be\la{0.4}
\int_{B_1} |\nabla u|^2 dx \geq \frac14   \int_{B_1} \frac{u^2}{(1-|x|)^2} dx + b_n  \int_{B_1} u^2 dx,
\ee
for some positive constant $b_n$. This time the best constant $b_n$  depends on the space dimension
with  $b_n> \mu_1$ when $n \geq 4$,
but   in the $n=3$     case, one has that $b_3 = \mu_1$, see \cite{BFT}.

On the  other  hand
the  classical  Sobolev  inequality
\be\la{0.5}
 \int_{\R^n} |\nabla u|^2 dx \geq   S_n
 \left( \int_{\R^n} |u|^{\frac{2n}{n-2}} dx \right)^{\frac{n-2}{n}},
\ee
is valid  for any  $u \in C^{\infty}_{0}(\R^n)$ where  $S_n=
 \pi n(n-2) \left( \Gamma(\frac{n}{2})/ \Gamma(n) \right)^{\frac2n}$ is the best constant, see \cite{A}, \cite{T}.
Maz'ya  \cite{M} combined  both the Hardy and Sobolev  term   in one inequality  valid
 in   the upper half space.  After   a conformal transformation it leads to the following
Hardy--Sobolev--Maz'ya  inequality
\be\la{0.6}
\int_{B_1} |\nabla u|^2 dx \geq \frac14   \int_{B_1} \frac{u^2}{(1-|x|)^2} dx + B_n \left( \int_{B_1}
 |u|^{\frac{2n}{n-2}} dx \right)^{\frac{n-2}{n}},
\ee valid for any  $u \in C^{\infty}_{0}(B_1)$.  Clearly  $B_n
\leq S_n$ and it was shown in \cite{TT} that $B_n < S_n$ when $n
\geq 4$.  Again,  the case   $n=3$   turns out to be special.
Benguria Frank and Loss  \cite{BFL}  have recently established
that $B_3=S_3=3 (\pi/2)^{4/3}$ (see also Mancini and Sandeep
\cite{MS}).

When distance is taken from the origin the analogue  of (\ref{0.6}) has been established
in \cite{FT} by methods quite different to the ones we use in the present work.
  To state the result we first define
 \be \la{1.log} X_1(a,s) :=( a- \ln s)^{-1}, \;\;\;\; a  > 0,~~~~ 0<s\leq1.
\ee
We then have:
\be\la{0.8}
\int_{B_1} |\nabla u|^2 dx \geq   \left(\frac{n-2}{2} \right)^2   \int_{B_1} \frac{u^2}{|x|^2} dx
 +
C_n(a) \left( \int_{B_1} X_1^{\frac{2(n-1)}{n-2}}(a,|x|) |u|^{\frac{2n}{n-2}}  dx
 \right)^{\frac{n-2}{n}}.
 \ee
We note that one cannot remove the logarithm  $X_1$ in (\ref{0.8}) and actually the exponent
 $\frac{2(n-1)}{n-2}$ is optimal.  Our main concern in this note is to calculate the
best constant $C_n(a)$ in (\ref{0.8}).  To this end we have:

\vspace{2mm}
\noindent {\bf Theorem A}{\em ~
Let  $n \geq 3$.
The best constant $C_n(a)$  in (\ref{0.8})  satisfies
\[
C_n(a)  = \left\{ \begin{array}{ll}
 (n-2)^{-\frac{2(n-1)}{n}}  \, S_{n}, &  ~~~~~~~~~  a \geq \frac{1}{n-2}  \\
 a^{\frac{2(n-1)}{n}}  \, S_{n}, &~~~~0< a < \frac{1}{n-2}.\\
\end{array} \right.
\]
  When restricted to radial functions,
the best constant in  (\ref{0.8}) is given by
\[
C_{n, radial} (a) = (n-2)^{-\frac{2(n-1)}{n}}  \, S_{n}, ~~~~~~{\rm for~all} ~~~a \geq 0.
\]
In all cases there is no $H_0^1(B_1)$ minimizer.
} \vspace{2mm}

\noindent
 One easily  checks that $C_n(a) < S_n$
when $n \geq 4$.  Surprisingly, in the $n=3$ case  one has that
$C_3(a)=S_3=3 (\pi/2)^{4/3}=B_3$, for $a \geq 1$, that is,
inequalities (\ref{0.5}), (\ref{0.6}) and (\ref{0.8}) share the
same best constant.

 Using the
change of variables  $u(x) = |x|^{-\frac{n-2}{2}} v(x)$
 inequality (\ref{0.8}) is easily seen to be equivalent to
\be
\int_{B_1} |x|^{-(n-2)}  |\nabla v|^2 dx
 \geq C_n(a)
\left( \int_{B_1} |x|^{-n}  X_1^{\frac{2(n-1)}{n-2}}(a,|x|) ~ |v|^{\frac{2n}{n-2}}  dx
 \right)^{\frac{n-2}{n}},  ~~~~~v \in C^{\infty}_{0}(B_1).
\label{1.15}
\ee
For later use
we denote by $\spw$ the completion of $ C^{\infty}_{0}(B_1)$  under the norm
 $\left( \int_{B_1} |x|^{-(n-2)}  |\nabla v|^2 dx  \right)^{1/2}$.

Estimate  (\ref{1.15}) is a limiting case of a Caffarelli--Kohn--Nirenberg inequality.
Indeed, for any  $-\frac{n-2}{2} < b < \infty$, the following  inequality holds:
\be\la{ckn}
\int_{\R^n} |x|^{2b}  |\nabla v|^2 dx   \geq  S(b,n)  \left(\int_{\R^n}
  |x|^{\frac{2bn}{n-2}} |v|^{\frac{2n}{n-2}} \right)^{\frac{n-2}{n}},  ~~~~~v \in C^{\infty}_{0}(\R^n);
\ee see \cite{CKN}, Catrina and Wang \cite{CW}. Moreover, for
$b=-\frac{n-2}{2}$ estimate (\ref{ckn})  fails.  Clearly, estimate
(\ref{1.15}) is the   limiting case of (\ref{ckn}) for
$b=-\frac{n-2}{2}$.  Thus we have:

\vspace{2mm}
\noindent {\bf Theorem A'}{\em ~
Let  $n \geq 3$.
The best constant $ C_n(a)$ in the limiting Caffarelli--Kohn--Nirenberg inequality (\ref{1.15})
is given
\[
C_n(a)  = \left\{ \begin{array}{ll}
 (n-2)^{-\frac{2(n-1)}{n}}  \, S_{n}, &  ~~~~~~~~~  a \geq \frac{1}{n-2}  \\
 a^{\frac{2(n-1)}{n}}  \, S_{n}, &~~~~0< a < \frac{1}{n-2}.\\
\end{array} \right.
\]
  When restricted to radial functions,
the best constant in  (\ref{1.15}) is given by
\[
C_{n, radial} (a) = (n-2)^{-\frac{2(n-1)}{n}}  \, S_{n}, ~~~~~~{\rm for~all} ~~~a \geq 0.
\]
In all cases there is no $\spw$ minimizer.
}
\vspace{2mm}

\noindent
We  note that the nonexistence of a  $\spw$ minimizer of Theorem A' is stronger than the nonexistence
of an  $H_0^1(B_1)$   minimizer of Theorem A.  This is due to the fact that the existence of an
 $H_0^1(B_1)$   minimizer  for (\ref{0.8})  would imply existence of a $\spw$ minimizer for (\ref{1.15}),
see  Lemma 2.1 of \cite{FT}.

The above results can be easily  transformed to  the exterior of the unit ball $B_1^c$. For instance we have:

\vspace{2mm}
\noindent {\bf Corollary}{\em ~
Let  $n \geq 3$. For any  $u \in C^{\infty}_{0}(B_1^c)$,   there  holds
\be\la{cor}
\int_{B_1^c} |\nabla u|^2 dx \geq   \left(\frac{n-2}{2} \right)^2   \int_{B_1^c} \frac{u^2}{|x|^2} dx
 +
C_n(a) \left( \int_{B_1^c} X_1^{\frac{2(n-1)}{n-2}} \left(a,\frac{1}{|x|} \right) |u|^{\frac{2n}{n-2}}  dx
 \right)^{\frac{n-2}{n}}.
 \ee
where the best constant  $C_n(a)$ is the same as in Theorem A.
}
\vspace{2mm}

Our method can also cover the case of a general  bounded domain $\Omega$ containing the origin. In particular
we have

 \vspace{2mm}
\noindent {\bf  Theorem  B}{\em~
Let  $n \geq 3$  and  $\Omega \subset  \R^n$ be  a bounded domain containing the origin.  Set $D := \sup_{x \in \Omega}|x|$.
 For any  $u \in C^{\infty}_{0}(\Omega)$, there holds
\be\la{thb}
\int_{\Omega} |\nabla u|^2 dx \geq   \left(\frac{n-2}{2} \right)^2   \int_{\Omega} \frac{u^2}{|x|^2} dx
 +
C_n(a) \left( \int_{\Omega} X_1^{\frac{2(n-1)}{n-2}} \left(a,\frac{|x|}{D} \right) |u|^{\frac{2n}{n-2}}  dx
 \right)^{\frac{n-2}{n}},
\ee
where the best constant $C_n(a)$ is independent of $\Omega$ and is  given by
 \[
C_n(a)  = \left\{ \begin{array}{ll}
 (n-2)^{-\frac{2(n-1)}{n}}  \, S_{n}, &  ~~~~~~~~~  a \geq \frac{1}{n-2}  \\
 a^{\frac{2(n-1)}{n}}  \, S_{n}, &~~~~0< a < \frac{1}{n-2}.\\
\end{array} \right.
\]
}

It follows easily from Theorem A' that there no minimizers for (\ref{cor}) and (\ref{thb}) in the
appropriate energetic  function  space.

 \vspace{2mm}
We next consider the $k$--improved  Hardy--Sobolev inequality derived in \cite{FT}.
Let  $k$ be a fixed  positive integer. For $X_1$ as in (\ref{1.log}) we define
 for $s \in (0,1)$,
\be\la{2.30}
X_{i+1}(a,s) = X_1(a, X_{i}(a,s)), \qquad i=1,2,\dots,k.
\ee
Noting that  $X_{i}(a,s)$ is a decreasing function of $a$ we easily check that
 there exist  unique positive constants
$0< a_k < \beta_{n,k} \leq 1$ such that :  \\
(i) The  $X_{i}(a_k,s)$  are well  defined for all  $i=1,2\ldots,k+1$, and all  $s \in (0,1)$
and  $X_{k+1}(a_k,1)=\infty$.  In other words, $a_k$ is the minimum  value of the  constant $a$
so that the  $X_{i}$'s, $i=1,2\ldots,k+1$,  are all well defined in $(0,1)$. \\
(ii) $X_1(\beta_{n,k}, 1) X_2(\beta_{n,k}, 1) \ldots  X_{k+1}(\beta_{n,k}, 1) = n-2$.

  For $n  \geq 3$, $k$ a fixed
positive integer and
 $u \in C^{\infty}_{0}(B_1)$ there holds:
 \bea
\la{0.11}
\int_{B_1}
|\nabla u|^2 dx \! \! & \geq             &   \left(\frac{n-2}{2} \right)^2
\int_{B_1} \frac{u^2}{|x|^2} dx
 + \frac14 \sum_{i=1}^{k}  \int_{B_1} \frac{X_1^2(a,|x|) \ldots X_{i}^2(a,|x|)}{|x|^2} u^2 dx   \nonumber \\
&   + & \! \! \!  C_{n,k}(a) \left( \int_{B_1} (X_1(a,|x|) \ldots
X_{k+1}(a,|x|))^{\frac{2(n-1)}{n-2}} |u|^{\frac{2n}{n-2}}  dx
 \right)^{\frac{n-2}{n}}.
\label{1.13}
\eea
In our next result we calculate the best constant $ C_{n,k}(a)$ in (\ref{0.11}).

\vspace{2mm}
\noindent {\bf Theorem C}{\em ~
Let  $n \geq 3$  and  $k=1,2,...$ be  a fixed  positive integer.  The best constant
 $ C_{n,k}(a)$ in (\ref{0.11})  satisfies:
\[
C_{n,k}(a)  = \left\{ \begin{array}{ll}
 (n-2)^{-\frac{2(n-1)}{n}}  \, S_{n}, &  ~~~~~~~~~~~  a \geq  \beta_{n,k}  \\
 \left( \prod_{i=1}^{k+1}X_{i}(a,1) \right)^{-\frac{2(n-1)}{n}}  \, S_{n}, &~~~~~~a_k< a <\beta_{n,k}.\\
\end{array} \right.
\]
  When restricted to radial functions,
the best constant of (\ref{1.13}) is given by
\[
C_{n,k, radial} (a) = (n-2)^{-\frac{2(n-1)}{n}}  \, S_{n}, ~~~~~~{\rm for~all} ~~~ a > a_k.
\]
}
 \vspace{2mm}

\noindent
Again we notice that $C_{n,k}(a) < S_n$ for  $n \geq 4$  but  $C_{3,k}=S_3$ for
$a \geq  \beta_{3,k}$.

As in Theorem A, one can establish  by similar arguments
 the  nonexistence of an  $H_0^1(B_1)$ minimizer to (\ref{1.13}), as well as the  analogues
 of   Theorem A', Corollary  and Theorem B in the
case of the  $k$--improved  Hardy--Sobolev inequality.

\setcounter{equation}{0}
\section{The proofs}

Theorems A  follows from Theorem  A', we therefore prove  Theorem A':

\noindent
{\bf Proof of Theorem A':}  At first we will show that
\be\la{2.1}
C_{n} (a)  =  (n-2)^{-\frac{2(n-1)}{n}}  \, S_{n},
~~~~~~{\rm when}  ~~~ a \geq \frac{1}{n-2}.
\ee
We have that
\be
\la{2.2}
C_{n} (a)  = \inf_{v \in  C^{\infty}_{0}(B_1)} \frac{\int_{B_1} |x|^{-(n-2)}  |\nabla v|^2 dx}
{\left( \int_{B_1} |x|^{-n}  X_1^{\frac{2(n-1)}{n-2}}(a,|x|) ~ |v|^{\frac{2n}{n-2}}  dx
 \right)^{\frac{n-2}{n}}} \;.
\ee
We  change variables by  ($r=|x|$)
\be\la{2.3}
v(x) = y(\tau , \theta), \qquad \quad  \tau=\frac{1}{X_1(a,r)}= a - \ln r,
 ~~~~ \theta = \frac{x}{|x|}.
\ee
This change of variables maps the unit ball $B_1 =\{x: |x|<1\}$ to the
complement of the  ball of radius $a$, that is, $B_a^c= \{ (\tau, \theta): a < \tau < +\infty,   ~~\theta \in S^{n-1} \}$.
Noticing that $X_1'(a,r)=\frac{X_1^2(a,r)}{r}$, $d \tau = -\frac{X_1'(a,r)}{X_1^2(a,r)}= -\frac{dr}{r}$,
we also  have
\[
| \nabla v |^2 =
\left( \frac{\partial v}{\partial r} \right)^2 + \frac{1}{r^2} | \nabla_{\theta} v|^2 =
 e^{2(\tau-a)} ( y^2_\tau + | \nabla_{\theta} y|^2).
\]
A straightforward calculation shows that for $y \in C^{\infty}([a,
\infty) \times S^{n-1})$ under Dirichlet boundary condition on
$\tau=a$ we have \be\la{2.5} C_{n} (a)  = \inf_{y(a,\theta)=0}
\frac{ \int_a^{\infty} \int_{S^{n-1}}
 ( y^2_{\tau} + | \nabla_{\theta} y|^2) dS d\tau}
{ \left( \int_a^{\infty} \int_{S^{n-1}}
 \tau^{-\frac{2(n-1)}{n-2}} |y|^{\frac{2n}{n-2}}dS d\tau \right)^{\frac{n-2}{n}}}.
\ee

In the sequel we will  relate $C_{n} (a)$  with the best
Sobolev constant  $S_{n}$.  It is well known that
for any $R$ with  $0<R \leq \infty$,
 \be\la{2.7} S_{n} = \inf_{u
\in C_0^{\infty}(B_R)} \frac{ \int_{B_R} |\nabla u|^2 dx} {
\left(\int_{B_R} |u|^{\frac{2n}{n-2}} dx \right)^{\frac{n-2}{n}}}.
\ee We also know that $S_{n}=S_{n,radial}$ the latter being the
infimum when taken over radial functions.
Changing variables in (\ref{2.7}) by
\be\label{2.9}
u(x) = z(t, \theta), \qquad \quad t= |x|^{-(n-2)},  \qquad  \theta = \frac{x}{|x|},
\ee
it follows that for any $R\in (0, \infty]$,
\be\label{2.11}
 (n-2)^{-\frac{2(n-1)}{n}}  \, S_{n}
=      \inf_{z(R^{-(n-2)}, \theta)=0}
 \frac{
 \int_{R^{-(n-2)}}^{\infty} \int_{S^{n-1}}
 ( z^2_t +  \left( \frac{1}{n-2} \right)^2 \frac{1}{ t^2}
 | \nabla_{\theta} z|^2) dS dt
}
{
 \left( \int_{R^{-(n-2)}}^{\infty} \int_{S^{n-1}}
 t^{-\frac{2(n-1)}{n-2}} |z|^{\frac{2n}{n-2}} dS dt \right)^{\frac{n-2}{n}}
}.
\ee
We note that a function $u$ is radial in $x$ if and only if the function $z$ is a function
of $t$ only. Looking at   (\ref{2.5}) and (\ref{2.11})  we have that
\be\la{2.13}
C_{n} (a) \leq C_{n, radial} (a)  =  (n-2)^{-\frac{2(n-1)}{n}}  \, S_{n,radial} =
 (n-2)^{-\frac{2(n-1)}{n}}  \, S_{n}.
\ee
On the other hand  let   us  take  $R = a ^{-\frac{1}{n-2}}$  (so that $a= R^{-(n-2)}$)
and
 assume  that $a \geq \frac{1}{n-2}$. Then
 $\left( \frac{1}{n-2} \right)^2 \frac{1}{ t^2}  \leq  1$  since $t \geq a \geq \frac{1}{n-2}$,
and therefore
\[
 C_{n} (a)  \geq \left(  \frac{1}{n-2}  \right)^{\frac{2(n-1)}{n}} S_{n}.
\]
Combining this with (\ref{2.13}) we conclude our claim (\ref{2.1}).

Our next  step is to prove the following: For any $a>0$ we have that
\be\la{2.15}
  C_{n}(a)  \leq  a^{\frac{2(n-1)}{n}} \;   S_{n}.
\ee
To this end
let  $0 \neq x_0 \in B_1$   and consider  the sequence of functions
\be\la{2.20}
U_{\xe}(x) =  ( \xe + |x-x_0|^{2})^{-\frac{n-2}{2}}  \phi_{\delta}(|x-x_0|),
\ee
where   $\phi_{\delta}(t)$ is a  $C_0^{\infty}$ cutoff function  which is zero for
$t > \delta$ and equal to one for $t < \delta/2$;   $\delta$ is small enough
so that  $|x_0| +\delta <1$ and therefore $U_{\xe} \in  C_0^{\infty}(B_{\delta}(x_0))
\subset C_0^{\infty}(B_1)  $.

 Then,  it is well known,  cf  \cite{BN},  that
\be\la{2.40}
 S_{n} = \lim_{\xe \ra  0} \frac{ \int_{B_{1}}
|\nabla U_{\xe}|^2 dx} { \left(\int_{B_{1}}
|U_{\xe}|^{\frac{2n}{n-2}} dx  \right)^{\frac{n-2}{n}}}.
 \ee
 From (\ref{2.2})  we have  that  for any $\xe>0$ small enough,
 \bea
 C_{n}(a) &
=  &  \inf_{v \in C_0^{\infty}(B_1)} \frac{ \int_{B_1}  |x|^{-(n-2)}  |\nabla
v|^2 dx} { \left(\int_{B_1} |x|^{-n}  X_1^{\frac{2(n-1)}{n-2}}(a,|x|)
|v|^{\frac{2n}{n-2}} dx  \right)^{\frac{n-2}{n}}}
\nonumber \\
&  \leq  &
 \frac{ \int_{B_{\delta}(x_0)} |x|^{-(n-2)} |\nabla U_{\xe}|^2 dx}
{ \left(\int_{B_{\delta}(x_0)} |x|^{-n}  X_1^{\frac{2(n-1)}{n-2}}(a,|x|) |U_{\xe}|^{\frac{2n}{n-2}} dx  \right)^{\frac{n-2}{n}}}
\nonumber \\
&  \leq  &
\left( \frac{|x_0|+\delta}{|x_0|-\delta} \right)^{n-2}
\frac{1}{X_1^{\frac{2(n-1)}{n}}(a,|x_0| - \delta)} ~
\frac{ \int_{B_{\delta}(x_0)} |\nabla U_{\xe}|^2 dx}
{ \left(\int_{B_{\delta}(x_0)}  |U_{\xe}|^{\frac{2n}{n-2}} dx  \right)^{\frac{n-2}{n}}},
\nonumber
\eea
where we used the fact that $X_1(a,s)$ is an increasing function of $s$.
Taking the limit $\xe \ra 0$ we conclude:
\[
C_{n}(a)  \leq  \left( \frac{|x_0|+\delta}{|x_0|-\delta} \right)^{n-2}  \frac{S_{n}}{X_1^{\frac{2(n-1)}{n}}(a,|x_0| - \delta)}.
\]
This is true for any $\delta>0$ small enough, therefore
\[
C_{n}(a)  \leq   X_1^{-\frac{2(n-1)}{n}}(a,|x_0|) \;   S_{n}.
\]
Since  $|x_0|<1$   is   arbitrary  and $X_1(a,s)$ is an increasing function of $s$,
 we end up with
\be\la{2.42}
  C_{n}(a)  \leq   X_1^{-\frac{2(n-1)}{n}}(a,1) \;   S_{n} = a^{\frac{2(n-1)}{n}} \;   S_{n},
\ee
and this proves our claim (\ref{2.15}).

To complete the calculation of $C_{n}(a)$ we  will  finally  show that
\be
\la{2.43}
 C_{n}(a)  \geq  a^{\frac{2(n-1)}{n}} \;   S_{n},~~~~~~~{\rm when}  ~~~ 0< a <  \frac{1}{n-2}.
\ee
To prove this we will relate the infimum  $ C_{n}(a)$  to a
Caffarelli--Kohn--Nirenberg inequality.  We will  need the following  result:

\begin{prop}\la{lem2}
Let     $b>0$    and
\be\label{2.44}
S_{n}(b) := \inf_{v \in C_0^{\infty}(\R^n)}
 \frac{ \int_{\R^n}|x|^{2b} |\nabla u|^2 dx}
{ \left(\int_{\R^n} |x|^{\frac{2bn}{n-2}} |u|^{\frac{2n}{n-2}} dx  \right)^{\frac{n-2}{n}}}.
\ee
Then $S_{n}(b) = S_{n}$ and this constant is not achieved in the appropriate function space.
\end{prop}
This is proved in   Theorem 1.1  of  \cite{CW}.

We change   variables   in (\ref{2.44})  by
\be\label{2.46}
u(x) = z(t, \theta), \qquad \quad t= |x|^{-(n-2)-2b}, \qquad  \theta = \frac{x}{|x|}.
\ee
A straightforward calculation shows that for any $R'$,
\be\la{2.48}
\left(n-2+2b  \right)^{-\frac{2(n-1)}{n}}  S_{n}
 \leq   \inf_{z(R',\theta)=0} \frac{\int_{R'}^{\infty} \int_{S^{n-1}}
 \left( z^2_t +   \frac{1}{ \left(n-2+2b  \right)^2 t^2}  | \nabla_{\theta} z|^2 \right) dS dt}
{ \left(\int_{R'}^{\infty} \int_{S^{n-1}}
 t^{-\frac{2(n-1)}{n-2}} |z|^{\frac{2n}{n-2}} dS dt\right)^{\frac{n-2}{n}}} .
\ee
Taking $R'=a$ and
comparing  (\ref{2.48})   with (\ref{2.5}) we  have that if
\be\la{2.50}
1 \geq  \frac{1}{ \left(n-2+2b  \right)^2 t^2},~~~~~~~~ {\rm for} ~~~~ t \geq a,
\ee
then
\be\la{2.52}
C_{n}(a) \geq \left( n-2+2b  \right)^{-\frac{2(n-1)}{n}}  S_{n}.
\ee
Condition (\ref{2.50})  is satisfied  if we choose $b \in (0, +\infty)$ such that
\be\la{2.53}
 \frac{1}{n-2}  > a = \left( n-2+2b  \right)^{-1}>0.
\ee
For such a $b$ it follows from (\ref{2.52}) that
\[
C_{n}(a) \geq    a^{\frac{2(n-1)}{n}} \;   S_{n},
\]
and this proves our claim  (\ref{2.43}).

We finally establish the nonexistence of an  energetic  minimizer.  We will  argue by contradiction. Suppose that
$\bar{v} \in \spw$ is a minimizer  of (\ref{2.2}).
Through the change of variables (\ref{2.3}),  the quotient in (\ref{2.5}) admits also a minimizer  $\bar{y}$.

 Consider first the case when
$a  \geq \frac{1}{n-2}$.  Comparing  (\ref{2.5}) and (\ref{2.11})  with $R = a^{-\frac{1}{n-2}}$, we conclude
that $\bar{y}$ is a  radial  minimizer of  (\ref{2.11}) as well.  It then follows that (\ref{2.7}) admits
a radial $H_0^1(B_R)$ minimizer $\bar{u}(r) = \bar{y}(t)$,  $t= r^{-(n-2)}$,  which  contradicts
the fact that  the  Sobolev  inequality  (\ref{2.7})   has no $H_0^1$ minimizers.

In the case  when  $0< a  <\frac{1}{n-2}$, we use a similar argument comparing  (\ref{2.5})  and (\ref{2.48})
to conclude the existence of a  radial  minimizer to (\ref{2.48}) with $b$ as in (\ref{2.53}).
This contradicts the nonexistence of minimizer for (\ref{2.44}).
 The proof of Theorem A' is now complete.

\noindent
{\bf Proof of Corollary:} One can argue in a similar way as in the previous proof, or apply Kelvin
transform to the inequality of  Theorem A.

\noindent
{\bf Proof of Theorem B:}  The lower bound on the best constant follows from Theorem A, the fact that if
   $u \in C^{\infty}_{0}(\Omega)$ then   $u \in C^{\infty}_{0}(B_D)$  (since  $\Omega \subset B_D$)
 and a simple scaling argument.

 To establish
the upper bound in the case where $0<a <\frac{1}{n-2}$
 we argue exactly  as in the proof of (\ref{2.15}) using  the test functions (\ref{2.20})
that concentrate near a  point of the  boundary of $\Omega$, that realizes the $\max_{x \in \Omega} |x|$.
Let us now consider the case where  $a \geq  \frac{1}{n-2}$. For $a>0$ and $0<\rho<1$,  we set
\[
\tilde{C}_n(a,\rho):= \inf_{u \in C^{\infty}_{0}(B_{\rho})} \frac{\int_{B_{\rho}} |\nabla u|^2 dx -
 \left(\frac{n-2}{2} \right)^2   \int_{B_{\rho}} \frac{u^2}{|x|^2} dx}
{\left( \int_{B_{\rho}} X_1^{\frac{2(n-1)}{n-2}}(a,|x|) |u|^{\frac{2n}{n-2}}  dx
 \right)^{\frac{n-2}{n}}
} ~ .
\]
A simple scaling argument and Theorem A shows that:
\[
\tilde{C}_n(a, \rho)= C_n(a - \ln \rho).
\]
Thus, for $\rho$ small enough we have that
\[
\tilde{C}_n(a, \rho) =  (n-2)^{-\frac{2(n-1)}{n}}  \, S_{n}.
\]
Since for $\rho$ small,  $B_{\rho} \subset \Omega$ the upper bound follows
easily in this case as well.

\noindent
{\bf Proof of Theorem C:} To simplify the presentation we will  write $X_i(|x|)$
instead of $X_i(a,|x|)$.
Let  $k$ be a fixed  positive integer.  We first consider the case  $a \geq \beta_{k,n}$.
We  change variables in (\ref{1.13})  by
\[
u(x)= |x|^{-\frac{n-2}{2}} X_1^{-1/2}(|x|) X_2^{-1/2}(|x|) \dots X_k^{-1/2}(|x|) v(x),
\]
  to obtain
\[
\int_{B_1} |x|^{-(n-2)} X_1^{-1}(|x|)\dots X_k^{-1}(|x|) |\nabla v|^2 dx \geq
\]
\be
C_{n,k}(a)  \left( \int_{B_1}
 |x|^{-n} X_1(|x|)\dots  X_k(|x|)  X_{k+1}^{\frac{2(n-1)}{n-2}}(|x|)
 |v|^{\frac{2n}{n-2}}
\right)^{\frac{n-2}{n}}, ~~~~v \in C_0^{\infty}(B_1).
\label{4.1}
\ee
We further change variables by
\[
v(x) = y(\tau, \theta), \qquad \tau=\frac{1}{X_{k+1}(r)},  \qquad
  \theta = \frac{x}{|x|}~~~~~(r=|x|).
\]
This change of variables maps the unit ball $B_1 =\{x: |x|<1\}$ to the
complement of the  ball of radius $r_a:=X_{k+1}^{-1}(1)$, that is,
 $B_{r_a}^c= \{ (\tau, \theta):X_{k+1}^{-1}(1)  < \tau < +\infty,   ~~\theta \in S^{n-1} \}$.
Note that
\[
d \tau = -\frac{X'_{k+1}(r)}{X_{k+1}^2(r)} dr = - \frac{ X_1(r)\dots  X_k(r)}{r}dr.
\]
Let  us denote by  $f_1(t)$  the inverse function  of $X_1(t)$.
We also  set  $f_{i+1}(t) = f_1(f_i(t))$, $i=1,2,\dots,k$.
  Consequently, $r= f_{k+1}(\tau^{-1})$. Also, $X_1(r)=f_k(\tau^{-1})$,
 $X_2(r)=f_{k-1}(\tau^{-1}),$ $\ldots$  $X_k(r)=f_{1}(\tau^{-1})$.

We then find
\be\la{4.4}
C_{n,k}(a) = \inf_{y(r_a,\theta)=0} \frac{ \int_{r_a}^{\infty} \int_{S^{n-1}}
 ( y^2_{\tau} + \left(f_1(\tau^{-1}) \dots f_k(\tau^{-1}) \right)^{-2}
 | \nabla_{\theta} y|^2) dS d\tau}
{ \left( \int_{r_a}^{\infty} \int_{S^{n-1}}
 \tau^{-\frac{2(n-1)}{n-2}} |y|^{\frac{2n}{n-2}}dS d\tau \right)^{\frac{n-2}{n}}}.
\ee
Again, we will relate this with the best Sobolev constant $S_n$. From  (\ref{2.11})
we have that
\be\label{4.6}
(n-2)^{-\frac{2(n-1)}{n}} ~   S_{n} =      \inf_{z(r_a, \theta)=0}
 \frac{
 \int_{r_a}^{\infty} \int_{S^{n-1}}
 ( z^2_t +  \frac{1}{(n-2)^2  t^2}
 | \nabla_{\theta} z|^2) dS dt
}
{
 \left( \int_{r_a}^{\infty} \int_{S^{n-1}}
 t^{-\frac{2(n-1)}{n-2}} |z|^{\frac{2n}{n-2}} dS dt \right)^{\frac{n-2}{n}}
}.
\ee
Comparing this with (\ref{4.4}) we have that
\be\la{4.8}
C_{n,k}(a)  \leq C_{n,k,radial}(a) =(n-2)^{-\frac{2(n-1)}{n}}  S_{n,radial}
=(n-2)^{-\frac{2(n-1)}{n}} S_n.
\ee
On the other hand for $a \geq \beta_{k,n}$  and $\tau \geq r_a$ we have that
\bea
\left( \tau^{-1}   f_1(\tau^{-1}) \dots f_k(\tau^{-1})  \right)^{-2} & \geq &
\left( r_a^{-1} f_1(  r_{a}^{-1}) \dots f_k( r_{a}^{-1}) \right)^{-2}    \nonumber \\
 & =  & \left(X_1(a,1) \dots X_k(a,1) X_{k+1}(a,1) \right)^{-2}    \nonumber  \\
&  \geq &  \frac{1}{(n-2)^2},        \nonumber
 \eea
therefore
\[
\left( f_1(\tau^{-1}) \dots f_k(\tau^{-1})  \right)^{-2}
 \geq  \frac{1}{(n-2)^2  \tau^2}, ~~~~~~~~~~~\tau \geq r_a,
\]
and consequently,
\[
C_{n,k}(a)   \geq (n-2)^{-\frac{2(n-1)}{n}} S_n.
\]
>From this and (\ref{4.8}) it follows that
\[
C_{n,k}(a)  = (n-2)^{-\frac{2(n-1)}{n}} S_n,  ~~~~~{\rm when}~~~a \geq \beta_{k,n}.
\]

The case where $a_k <a <\beta_{k,n}$ is  quite similar to the case  $ 0< a <  \frac{1}{n-2}$
 in the proof of Theorem
A'. That is,  testing  in  (\ref{4.1}) the sequence   $U_{\xe}$ as defined in (\ref{2.20}),
  we first prove  that
\[
C_{n,k}(a)   \leq  \left( \prod_{i=1}^{k+1}X_{i}(a,1) \right)^{\frac{-2(n-1)}{n}}  \, S_{n},
\]
 by an argument quite similar to the one leading to  (\ref{2.42}).
Finally,  in the case $a_k <a <\beta_{k,n}$, we
 obtain the opposite inequality by
  comparing  the infimum in (\ref{4.4}) with the   infimum
in (\ref{2.48}). This time we take $R'=r_a$  and $b>0$ is chosen so that
\[
  \prod_{i=1}^{k+1}X_{i}(a,1) =n-2+2b.
\]
We omit further details.

\medskip
\noindent {\bf Acknowledgments} Adimurthi is thanking the
Departments of Mathematics and Applied Mathematics of University
of Crete for the invitation as well as the warm hospitality. The
authors thank the referee for raising the question of existence or
nonexistence of minimizers.

 \end{document}